\documentclass[a4paper,twoside]{amsart}
\usepackage[T1]{fontenc}

\newtheorem{theorem}{Theorem}[section]
\newtheorem{lemma}[theorem]{Lemma}
\newtheorem{proposition}[theorem]{Proposition}

\theoremstyle{definition}
\newtheorem{definition}[theorem]{Definition}

\theoremstyle{remark}
\newtheorem{remark}[theorem]{Remark}

\usepackage{amsmath}
\usepackage{amsfonts}
\usepackage{amssymb}

\newcommand{\bbC}{\mathbb{C}}
\newcommand{\bbR}{\mathbb{R}}

\newcommand{\dr}{\mathrm{d}}

\numberwithin{figure}{section}

\newcommand{\mon}{\mathcal{M}\mathrm{on}}

\newcommand{\cP}{\mathcal{P}}
\newcommand{\cA}{\mathcal{A}}

\newcommand{\ve}{\varepsilon}
\newcommand{\tht}{\theta}
\newcommand{\id}{\mathrm{I}}
\newcommand{\Ad}{\mathrm{Ad}}
\title[Contour integrals as $\Ad$-invariant functions]{Contour integrals as $\Ad$-invariant functions on the fundamental group}
\author{Marcin Bobie\'nski} 
\address{Institute of Mathematics, Warsaw University, ul. Banacha 2, 02-097 Warsaw, Poland.}
\email{mbobi@mimuw.edu.pl} 
\thanks{This research was supported by the KBN Grant No 2 P03A 010 22}

\begin{document}
\begin{abstract}
We introduce a general approach to contour integrals. It covers usual Abelian integrals, the higher order Melnikov integrals and the generalized Abelian integrals (see \cite{bozoell,bozo2}). We prove that the generating function always satisfies a linear differential equation of finite order. We also present a relationship between the generalized Abelian integral and certain representation of the fundamental group of complex curve.
\end{abstract}

\maketitle

\section{Introduction}
\label{sec:intro}

We start with recalling known examples of contour integrals which appear naturally in study of limit cycles of polynomial vector fields.

Analyzing the phase portrait of the planar polynomial perturbation of Hamiltonian vector field
\begin{equation}
  \label{2dsys}
  \begin{split}
    \dot{x} &= H_{y}+\ve P,\\
    \dot{y} &= -H_{x}+\ve Q
  \end{split}
\end{equation}
the notion of the Abelian integral
\[
I_\gamma = \int_\gamma Q \dr x - P \dr y
\]
naturally appears. It is a linear (in $\ve$) approximation of the Poincar\'{e} return map along an oval $\gamma\subset \{H= \mathrm{const}\}\subset\bbR^2$. The zeroes of Abelian integrals corresponds to the limit cycles of (\ref{2dsys}), providing $I_\gamma$ does not vanish identically. If $I_\gamma\equiv 0$, we study the higher order terms (with respect to the perturbation parameter $\ve$) in the expansion of the Poincar\'{e} return map (see \cite{zolq}). The first non-vanishing term $M_\gamma$ is called \emph{generating function} \cite{gi} (or \emph{higher order Melnikov integral}) and it detects the limit cycles of (\ref{2dsys}). 

The last contour integral, which we recall, is the \emph{generalized Abelian integral} introduced and investigated in \cite{lezo,bozoell,bozo2}. It appears naturally in study of polynomial perturbation of a multi-dimensional system with invariant plane on which it is Hamiltonian (see (\ref{mdsys}) below). The invariant surface becomes slightly deformed but survives, assuming certain normal hyperbolicity condition. The Poincar\'e return map of the system restricted to this surface is considered. The linear term of this map turns out to be the sum of the usual Abelian integral and certain new contour integral $\Psi_\gamma$ called \emph{generalized Abelian integral} (see formula (\ref{Psidef}) below).

Estimating the number of zeroes of these contour integrals is the important task in study of the infinitesimal Hilbert's 16th problem and its generalizations. This would give an upper bound for the number of limit cycles bifurcating from ovals of the Hamiltonian vector field. There are several approaches to the problem; in some of them crucial role plays the analytical continuation of the contour integral to a complex plane, possibly with some points removed. The main tool there is the analysis of the monodromy of this continuation. The source of the non-trivial monodromy of the contour integral is the monodromy of the oval $\gamma$ prolongated to the complex domain. More precisely, let $\Delta\subset \bbC$ be the set of complex critical values of the Hamiltonian $H$. The continuation $\gamma_h$ of the real oval $\gamma_h\subset H^{-1}(h)$ as $h$ varies along a curve in $\bbC\setminus\Delta$ is unique up to the homotopy. Thus the possible monodromy appears after surrounding the critical values from $\Delta$.

To describe the monodromy of contour integral around a critical value we need to know how does the loop $\gamma_h$ change and how does the integral depend on the loop. The change of the contour is described by the classical homotopical Picard-Lefshetz formula. The goal of this paper is to investigate dependence of the mentioned above contour integrals on the loop $\gamma$.

It is well known that the Abelian integral $I_\gamma$ is an additive function of $\gamma$. So, the map
\[
\pi_1(H^{-1}(h),*)\ni\gamma \mapsto I_\gamma
\]
factorizes through the homology group of the complex level curve $H_1(H^{-1}(h))$ and the monodromy of the Abelian integral describes (and defines) the Gauss-Manin connection on the Milnor homology bundle associated to the polynomial fibration 
\[
H: \bbC^2\setminus H^{-1}(\Delta) \longrightarrow \bbC\setminus\Delta.
\]
Moreover, since the homology group is finitely generated, the function $I_\gamma$ satisfies a differential equation of finite order.

The $k$-th order Melnikov integral (or the generating function) $M^k_\gamma$ is also an additive function of the contour $\gamma$, which belongs to certain subgroup $S$ of the fundamental group (see \cite{gi}). This subgroup consist of those loops for which all lower order terms in the expansion of the Poincar\'e return map vanish. We prove that the homomorphism $\gamma\mapsto M^k_\gamma$, $\gamma\in S$ factorizes through some finitely generated abelian group  (see Proposition \ref{pr:fingen} below). So, the generating function also satisfies a linear differential equation of finite order (Theorem \ref{th:de}). In general, the order of this equation depends on the system (\ref{2dsys}). 

In the last example of the contour integral we introduce a relationship between the generalized Abelian integral $\Psi$ and certain representation of the fundamental group of complex level curve $H^{-1}(h)$ (in Theorem \ref{th:rep}). This description results in explanations of the monodromy of generalized Abelian integrals. As a further application we prove (in Proposition \ref{pr:ellneq}) that the generalized Abelian integral does not satisfy any linear differential equation of finite order. Thus the generalized Abelian integrals are essentially different from the Abelian integrals and higher order Melnikov functions.

\section{The general contour integral}
\label{sec:genint}
\newcommand{\crc}{B}
\newcommand{\cpr}{\tau}

Let $H(x,y)$ be a polynomial (not necessarily real) with isolated critical points. We denote by $\Delta\subset \bbC$ the set of complex critical values of $H$ and 
\begin{gather*}
\crc \colon = \bbC\setminus \Delta, \qquad E\colon = \bbC^2\setminus H^{-1}(\Delta), \\
E_h\colon = \{(x,y)\in\bbC^2 : H(x,y)=h\} \subset E, \quad h\in \crc.
\end{gather*}

Let us recall the terminology related to the classes of loops on a manifold defined modulo homotopy; we shall use them below many times. Let $*\in E_0$ be a base point chosen on a connected manifold $E_0$. We denote by $\pi_1(E_0,*)$ the \emph{fundamental group} i.e.\ the group of loops starting at $*$ defined modulo homotopy with the base point fixed. The set of loops on $E_0$ defined modulo \emph{free} homotopy is denoted by $\pi_1(E_0)$ (this set does not have a group structure). We have the following canonical map (surjection)
\begin{equation}
  \label{canpr}
  \cpr: \pi_1(E_0,*) \rightarrow \pi_1(E_0).
\end{equation}
Note that two elements $\gamma,\delta\in \pi_1(E_0,*)$ are unified under the projection $\cpr$ iff they are related by an adjoint transformation
\[
\cpr(\gamma)=\cpr(\delta) \quad \Leftrightarrow \quad \delta = \Ad_\theta\, \delta = \theta\cdot\gamma\cdot \theta^{-1},
\]
with an element $\theta\in \pi_1(E_0,*)$.
\smallskip

The polynomial fibration $H: E\rightarrow \crc$ is locally trivial and, by the covering homotopy theorem (see\ \cite{ffg}), a closed loop $\gamma\subset E_{h_0}$ may be prolongated along any curve $l:[0,1]\rightarrow \crc$, $l(0)=h_0$, in the base to a 1-parameter family $\gamma_t\subset E_{l(t)}$ of loops; each loop $\gamma_t$ is determined uniquely up to homotopy. So any polynomial (as above) defines a covering
\[
\pi_1(E) \longrightarrow \crc,
\]
whose fibre over $h\in\crc$ is the set $\pi_1(E_h)$ of \emph{free} homotopy classes of loops on the level curve $E_h$ (see \cite{gi} for more details). This covering has many connected components, e.g.\ there always exists a component (isomorphic to $\crc$), corresponding to the class of contractible loops.

\begin{definition}
\label{df:ci}
A holomorphic (respectively meromorphic) \emph{contour integral} $K$ related to the polynomial fibration $H$ is a holomorphic (resp.\ meromorphic) function on some connected components of $\pi_1(E)$. The respective multi-valued function on $\crc$ will be denoted by $K_\gamma$, where $\gamma\subset E_{h_0}$ determines the connected component of $\pi_1(E)$.
\end{definition}

In important examples the contour integrals are defined by integrals along pa\-ra\-met\-rized loops. Let us present local properties of the contour integrals. Let $U\subset \crc$ be a contractible subset. We choose a base point $h_0\in U$ and $*\in E_0\colon = E_{h_0}$. Since the prolongation of a loop along $U$ is unique up to the homotopy, the covering $\pi_1(E)$ restricted to $U$ trivializes to the product
\[
\pi_1(E)|_{U} \cong \pi_1(E_0) \times U.
\]
Note that the fundamental group of the base $\pi_1(\crc,h_0)$ acts naturally on the fibre $\pi_1(E_0)$ via the monodromy.

The local description of a contour integral is given in the following.
\begin{proposition}
\label{pr:ci}
Let $K$ be a holomorphic (resp. meromorphic) contour integral related to the polynomial fibration $H$ and let $U\subset\crc$ be a contractible neighborhood of $h_0$. Then the following objects uniquely determine local properties of the contour integral $K$:
\begin{enumerate}
  \item A subset $S\subset \pi_1(E_0,*)$ which is $\Ad$-invariant and the image $\tau(S)\subset \pi_1(E_0)$ is monodromy invariant.
  \item An $\Ad$-invariant map 
    \begin{equation*}
      \kappa : S \rightarrow \cA(U), \qquad \gamma\mapsto \kappa(\gamma)=K_\gamma,
    \end{equation*}
    where $\cA(U)$ is the space of holomorphic (resp.\ meromorphic) functions on $U$. 
\end{enumerate}
\end{proposition}
\begin{proof} 
The contour integral $K$ is defined on some components of the covering $\pi_1(E)$. Locally, over the contractible neighborhood $U$ of $h_0$, the domain of definition $K$ has the form $\widetilde{S}\times U$, where $\widetilde{S}\subset \pi_1(E_0)$. This subset $\widetilde{S}$ is monodromy invariant, since we deal with a whole connected components of the covering $\pi_1(E)\rightarrow \crc$.

The contour integral $K$ defines the map $\widetilde{\kappa}:\widetilde{S} \rightarrow \cA(U)$. Taking the pre-image of $\widetilde{S}$ under the canonical map (\ref{canpr}), we get the subset $S\subset \pi_1(E_0,*)$. Composing $\widetilde{\kappa}$ with the map (\ref{canpr}) we get the map $\kappa$.

One can reconstruct the contour integral $K$ on $U$ from $\kappa$, using the relation $\kappa(\gamma) = K_\gamma$.\\
\end{proof}

\begin{remark}
\label{rk:mon}
In many interesting examples the function $\kappa$ is defined on a proper subset $S$ of the fundamental group. In practice, the subset $S$ is defined via another contour integral $K^0$. Locally on $U$ we have a pair $(\kappa^0,S^0)$ as described in Proposition \ref{pr:ci}. The subset $S=\{\gamma\in S^0:\ \kappa^0(\gamma)\equiv c\}, c\in\bbC$, is $\Ad$-invariant due to the $\Ad$-invariance of the function $\kappa^0$ and it is \emph{monodromy invariant} due to the unique continuation of analytical functions. This subset $S$ is the domain of definition of another contour integral. One can start this inductive procedure from the Abelian integral, which is defined for arbitrary loop and so the corresponding function $\kappa_I$ is defined on whole fundamental group $\pi_1(E_0,*)$.

The above inductive algorithm is in the spirit of J. P. Fran{\c{c}}oise procedure \cite{jpfr,jpfr2}.
\end{remark}

\subsection{The Abelian integrals}
\label{sec:ab}

Let $H$ be a polynomial with isolated critical points and let $\omega$ be a polynomial 1-form
\begin{equation}
  \label{omega}
  \omega = Q(x,y) \dr x - P(x,y) \dr y, \qquad P,Q\in\bbC[x,y].
\end{equation}
The Abelian integral 
\begin{equation}
  \label{abint}
  I_\gamma (h) = \int_\gamma \omega, \qquad \gamma\subset E_h,
\end{equation}
restricted to the neighborhood $U\ni h_0$, defines the additive function $\kappa_I$ on $S=\pi_1(E_0,*)$. So, this function factorizes through the homology group
\[
\kappa_I : \pi_1(E_0,*) \rightarrow H_1(E_0) \rightarrow \cA(U).
\]

\section{The higher order Melnikov integrals}
\label{sec:genfn}

We fix, as previously, a polynomial $H(x,y)$ with isolated critical points and a polynomial 1-form $\omega$ (\ref{omega}). We consider a holomorphic foliation defined by
\[
\dr H - \ve \omega =0.
\]
To any closed loop $\gamma\subset E_0$ we can assign a germ of holomorphic map $\cP_\ve (h)$ defined on a neighborhood of $h_0\in\crc$. If $\gamma$ is a real oval of the real polynomial $H$, then the map $\cP$ is just the complexification of the Poincar\'{e} return map. For the details of definition and proofs of some properties in general case we refer the reader to \cite{gi}. The map $\cP_\ve$ has the following properties.
\begin{enumerate}
\item The expansion of $\cP_\ve$ in $\ve$ has the form
\[
\cP_\ve (h)= h + \ve^k M^k_\gamma (h) +\ldots.
\]
In general, the index $k$ depends on the loop $\gamma$.
\item The \emph{generating function} $M^k_\gamma$ is the holomorphic contour integral related to the fibration $H$.
\item If $k=1$ then the generating function is equal to the Abelian integral (\ref{abint})
\[
M^1_\gamma = I_\gamma.
\]
\item Let $\gamma,\delta\subset E_0$ are two loops with common base point $*$ and the same index of expansion $k$. We denote by $\gamma\cdot \delta$ the composition of $\gamma$ and $\delta$ in $\pi_1(E_0,*)$ i.e.\ the loop consisting of $\gamma$ followed be $\delta$. The generating function of the composition is
  \begin{equation}
    \label{genad}
    M^k_{\gamma\cdot\delta} = M^k_{\gamma} + M^k_{\delta}.
  \end{equation}
\end{enumerate}

We prove the following (compare Theorem 2 in \cite{gi}).
\begin{theorem}~
\label{th:de}
Any generating function $M^k_\gamma$ satisfies a linear differential equation
\begin{equation}
  \label{gende}
  a_n x^{(n)} + a_{n-1} x^{(n-1)}+ \cdots + a_0 x =0,
\end{equation}
where $a_j$ are analytic functions on $\crc$.
\end{theorem}
\begin{remark}
  In general, the order of equation (\ref{gende}) depends on $k$. 
\end{remark}
\begin{remark}
  If, in addition to the assumptions of Theorem \ref{th:de}, one imposes the regular (i.e\ of power type) growth of the generating functions $M^k$, then the equation (\ref{gende}) is of the Fuchs type.
\end{remark}

\noindent
\begin{proof}[Proof of Theorem \ref{th:de}]
We fix a contractible neighborhood $U \subset \crc$. To simplify the notation in the subsequent proof we shall denote the fundamental group by
\[
\pi_1 \colon = \pi_1(E_0,*), \qquad E_0 = E_{h_0}, \quad h_0\in U.
\]

According to the Proposition \ref{pr:ci}, the generating function $M^k$ on $U$ of order equal to $k$ corresponds to the $\Ad$-invariant function
\[
m^k : S^k \rightarrow \cA(U), \quad \gamma\mapsto m^k(\gamma)=M^k_\gamma.
\]
It turns out that the subset $S^k$ is a subgroup of $\pi_1$.
We denote by $[\pi_1,S^k]$ the commutant group generated in $\pi_1$ by elements $\{g s g^{-1} s^{-1},\ g\in \pi_1,\ s\in S^k\}$. 

The proof of Theorem \ref{th:de} is based on the following proposition, which will be proved later.
\begin{proposition}
\label{pr:fingen}
The domain $S^k$ is a subgroup of $\pi_1$ and the quotient
\[
S^k/[\pi_1,S^k]
\]
is a finitely generated abelian group.
\end{proposition}
Now we finish the proof of Theorem. The function $m^k : S^k \rightarrow \cA(U)$ is $\Ad$-invariant (by definition) and is additive (by \ref{genad}). So, the map $m^k$ factorizes through the quotient group
\[
m^k : S^k \rightarrow S^k/[\pi_1,S^k] \rightarrow \cA(U).
\]
Since this quotient is finitely generated (by Proposition \ref{pr:fingen}), the contour integral $M^k$ belongs to some finite dimensional monodromy representation. Let us denote by $M^k_{\gamma_1}, \ldots ,M^k_{\gamma_n}$ the basis of the space $m^k(S^k)$. Any $M^k_{\gamma_j}$ is the contour integral and so it is defined globally as a multi-valued function on $\crc$. The vector space $\mathrm{Span}_{\bbC}(M^k_{\gamma_1}, \ldots ,M^k_{\gamma_n})$ is monodromy invariant due to the invariance of $S^k$. 

We consider the following equation
\[
\det\begin{pmatrix}
M^k_{\gamma_1} & (M^k_{\gamma_1})^\prime & \dots & (M^k_{\gamma_1})^{(n)}\\
M^k_{\gamma_2} & (M^k_{\gamma_2})^\prime & \dots & (M^k_{\gamma_2})^{(n)}\\
\hdotsfor{4} \\
M^k_{\gamma_n} & (M^k_{\gamma_n})^\prime & \dots & (M^k_{\gamma_n})^{(n)}\\
x & x^\prime & \dots & x^{(n)}\\
\end{pmatrix} =0.
\]
Taking the expansion of the determinant with respect to the last row, we get the equation of the form (\ref{gende}) with coefficients $a_j$ being determinants of suitable minors. The monodromy along any loop in $\pi_1(\crc,h_0)$ results in multiplication of all the coefficients by a common factor. Thus, any quotient $a_j/a_k$ is monodromy invariant and so meromorphic. Moreover, since the contour integrals $M^k_{\gamma_j}$ are defined globally, the coefficients $a_j$ of the equation (\ref{gende}) are meromorphic functions on the whole $\crc$. Multiplying the equation by a suitable holomorphic function, one can remove the poles of the coefficients.\\
\end{proof}

\noindent
\begin{proof}[Proof of Proposition \ref{pr:fingen}]
Let us introduce the notations:
\begin{align*}
  \pi_1^{(0)}   &= \pi_1= \pi_1(E_0,*),\\
  \pi_1^{(n+1)} &= [\pi_1,\pi^{(n)}].
\end{align*}
We show, by induction with respect to $k$, that $S^k$ is a subgroup of $\pi_1$ containing $\pi_1^{(k-1)}$. 

For $k=1$ the domain is the whole fundamental group $S^1 = \pi_1 = \pi_1^{(0)}$. 

The domain of the contour integral for $k+1$ is the kernel of $m^k$
\[
S^{k+1} = \ker (m^k) = \{\gamma \in S^k :\ m^k(\gamma) \equiv 0\}.
\]
By (\ref{genad}), the map $m^k$ is an additive map and so the kernel is a subgroup. Since $m^k$ is $\Ad$-invariant, we have the implication
\[
\pi_1^{(k-1)}\subset S^k \quad \Rightarrow \quad \pi_1^{(k)}\subset \ker(m^k) = S^{k+1}.
\]

Now we show that the quotient $S^k/[\pi_1,S^k]$ is finitely generated. Indeed, the quotient $\pi_1/\pi_1^{(k)}$ is a finitely generated nilpotent group (called super-solvable, see chapter 10 in \cite{hall}). This property is inherited by subgroups (see Theorem 10.5.1 in \cite{hall}). Since $\pi_1^{(k-1)}\subset S^{(k)}$, we have
\[
S^k/[\pi_1,S^k] = G/H,\qquad G= S^k/\pi_1^{(k)},\quad H=[\pi_1,S^k]/\pi_1^{(k)}.
\]
The group $G$ is finitely generated by the cited theorem, so the same is true for the quotient $S^k/[\pi_1,S^k]$. This quotient is also an abelian group since $[S^k,S^k]\subset [\pi_1,S^k]$. This finishes the proof of the proposition.\\
\end{proof}

\begin{remark}
In the above proof of Proposition \ref{pr:fingen} the crucial role is played by the Abelian group $S^k/[\pi_1,S^k]$. It depends, in general, on the Hamiltonian $H$, on $k$ and on the polynomial 1-form $\omega$ given by (\ref{omega}). This group is common for all loops $\gamma$ for which the leading term in the expansion of the Poincar\'{e} return map has order equal to $k$. 

On the other hand, L. Gavrilov and I. D. Iliev introduced in \cite{gi} another abelian group describing the monodromy of the generating function. It depends on the Hamiltonian $H$ and on the loop $\gamma$. Let us briefly recall their construction. We consider the loop $\gamma$ and all its prolongations along loops on the base $\crc$, i.e.\ the orbit of $\gamma$ under the monodromy action of the fundamental group of the base $\pi_1(\crc,h_0)$. Then we take the subgroup $\hat{S}$ generated by the pre-image of this orbit in the fundamental group $\pi_1(E_h,*)$. The following group
\[
H_\gamma \colon = \hat{S}/[\pi_1(E_h,*),\hat{S}]
\]
is abelian and any generating function gives rise to the homomorphism
\[
M: H_\gamma \rightarrow \cA(U),
\]
where $\cA(U)$ is the space of holomorphic functions on certain neighborhood $U$. The above homomorphism factorizes through the quotient $S^k/[\pi_1,S^k]$, so there exist a homomorphism
\[
H_\gamma \rightarrow S^k/[\pi_1,S^k]
\]
which is neither injective nor surjective in general. So, the method used in the above proof of Proposition \ref{pr:fingen} does not solve the open problem posted in \cite{gi}, whether the abelian group $H_\gamma$ is always finitely generated.
\end{remark}

\begin{remark}
In the most recent paper \cite{newgav} L. Gavrilov also has proved that the generating function always satisfies a linear differential equation of Fuchs type.
\end{remark}

\section{Generalized Abelian integrals}
\label{sec:nonab}

In all previous examples local properties of the contour integrals were described by the additive homomorphism defined on some subgroups of the fundamental group $\pi_1(E_0,*)$. Now we consider the generalized Abelian integrals which appear naturally in the study of certain multi-dimensional dynamical systems (see \cite{lezo,bozoell,bozo2}). Below we show that the generalized Abelian integral is an example of the contour integral whose local properties are described by a \emph{non-linear}, $\Ad$-invariant meromorphic function on certain subset (not a subgroup) of the fundamental group $\pi_1(E_0,*)$.
\medskip

We consider the following dynamical system which describes the perturbation of the system with invariant plane ($\{z=0\}$) on which it is Hamiltonian:
\begin{equation}
\label{mdsys} 
\begin{split}
\dot{x} &= H_{y}+z R+\ve P, \\
\dot{y} &= -H_{x}+zS+\ve Q, \\
\dot{z} &= A z+\ve b.
\end{split}
\end{equation}
All the functions $H,P,Q,R,S,A,b\in\bbR[x,y]$ are real polynomials and the Hamiltonian $H$ admits the region $D\subset\bbR^2=\{z=0\}$ filled with closed ovals of $H$. Imposing the \emph{normal hyperbolicity} condition on $D$ guarantees that the invariant surface $\{z=0\}$ deforms itself, but it survives the perturbation. 
The perturbed invariant surface $L_\ve$ has the following expansion in $\ve$
\[
L_\ve =\{(x,y,z): \ z= \ve g(x,y)+ \cdots ,\ (x,y)\in D\},
\]
where $g$ is the unique periodic (on ovals $\{H=h\}$ in $D$) solution to the \emph{normal variation equation}
\[
H_y \frac {\partial g}{\partial x} - H_x \frac{\partial g}{\partial y} = A g + b.
\]
An oval $\gamma$ has a natural parametrization by the Hamiltonian time $t$, given by the integral along $\gamma$ of the following 1-form
\begin{equation}
  \label{dt}
  \dr t = \left[\frac{\dr x}{H_y}\right]=\left[-\frac{\dr y}{H_x}\right] \mod \dr H.
\end{equation}
The period of $t$ on $\gamma$ is denoted by $T_\gamma=\int_\gamma \dr t$. The normal variation equation restricted to $\gamma$ has the form $\dot{g}(t) = A g + b$ and the periodic solution is given by the formula (see \cite{bozo2})
\begin{equation}
  \label{gintform}
  g(t) = \Big(e^{-\int_0^{T_\gamma} A(s)\dr s}-1\Big)^{-1} \int_t^{t+T_\gamma} b(s)\cdot e^{-\int_t^s A(u)\dr u}\,\dr s.
\end{equation}

Now we consider the Poincar\'{e} return map of the system (\ref{mdsys}) and its restriction to the invariant surface $L_\ve$. The Hamiltonian $H$ defines a coordinate on a transversal to the ovals. With respect to this coordinate, the return map has the following expansion in $\ve$
\begin{equation}
\label{pexp}
  \Delta H (h) = \ve J(h) + \cdots,\qquad J(h) = \int_{\gamma_h}(Q+g\cdot S)\dr x-(P+g\cdot R)\dr y,
\end{equation}
where $J(h)$ is the \emph{Pontryagin-Melnikov integral} (see \cite{pont,mel}). The integral $J(h)$ is given by a sum of the usual Abelian integral (\ref{abint}) and a new contour integral $\Psi_\gamma$, which has the form
\begin{equation*}
  \Psi_\gamma (h) = \int_{\gamma_h} (g \cdot p)\, \dr t,
\end{equation*}
where $p$ is a polynomial $p=(S\, H_y+R\, H_x)$ and $\dr t$ is given by (\ref{dt}).
\medskip

The above consideration motivates the following definition of \emph{generalized Abelian integral}. Let $H$ be a polynomial with isolated critical points and let $A,b,p\in\bbC[x,y]$ be three polynomials. For a contour $\gamma$ contained in the fibre $E_h$ of the polynomial fibration $H$, we define
\begin{equation}
\label{igamdef}
  I_\gamma \colon = \int_\gamma A\, \dr t, \qquad T_\gamma \colon = \int_\gamma \dr t.
\end{equation}
The integral of $\dr t$ is the multi-valued complex function $\gamma \overset{t}{\rightarrow} \bbC$. It can be extended to a holomorphic function on the neighborhood of $\gamma$ in $\bbC^2$.
\begin{lemma}
\label{lem:tholo}
For any loop $\gamma$ contained in the fibre $E_0$ there exist an open neighborhood $\gamma\subset V\subset \bbC^2$ and a multi-valued holomorphic function 
\[
t:\widetilde{V} \rightarrow \bbC
\]
such that the restriction of $\dr t$ to fibers coincides with the the Hamiltonian time 1-form (\ref{dt}) i.e.
\begin{equation}
  \label{dtv}
  H^\prime_y t^\prime_x - H^\prime_x t^\prime_y = 1.
\end{equation}
The increment of the coordinate $t$ along loop $\gamma$ is equal to the period $T_\gamma$.
\end{lemma}
\begin{proof}
We construct the function $t$ directly as integrals along paths in the fibers $E_h$. 

Let us choose a base point $*\in E_0$ and a small disk $D$ passing through $*$ and transversal to the fibers $E_h$. Any path $l$ contained in the fibre and starting on $D$
\[
l:[0,1]\rightarrow E_h, \quad l(0)\in D,
\]
provides a number $\int_l \dr t$, where $\dr t$ is the Hamiltonian time 1-form (\ref{dt}). The form $\dr t$ is holomorphic in sufficiently small neighborhood of $\gamma$, so the value of the integral $\int_l \dr t$ depends on the end points $l(0),\,l(1)$ and the homotopy class of $l$. The dependence on $(l(0),l(1))$ is \emph{holomorphic}. Any point $p$ sufficiently close to $\gamma$ belongs to the fibre $E_h$ which crosses $D$ in unique point $p_0$. Joining $p_0$ and $p$ by a path $l_p$ we get a function
\[
t(p)\colon = \int_{l_p} \dr t,
\]
which depends on the homotopy class of chosen $l_p$. Thus it is a multi-valued holomorphic function on a tubular neighborhood of $\gamma$. On a fibre $E_h$ the differential of $t$ coincide with (\ref{dt}). Since the Hamiltonian field $X_H=H^\prime_y \partial_x - H^\prime_x \partial_y$ is tangent to the fibre and satisfies $<\dr t,X_H>=1$, the identity (\ref{dtv}) follows.

The last statement about the period is obvious.\\
\end{proof}

We consider the following integral
\begin{equation}
  \label{Psidef}
  \Psi_\gamma (h) \colon = \Big(e^{-I_\gamma}-1\Big)^{-1} \int_0^{T_\gamma}\dr t  \int_t^{t+T_\gamma}\dr s\; p(t)\cdot b(s)\cdot e^{-\int_t^s A}.
\end{equation}

\begin{proposition}
\label{pr:gabi}
The formula (\ref{Psidef}) defines a meromorphic contour integral $\Psi$. The domain of definition of $\Psi$ consists of those components of $\pi_1(E)$ for which
\[
e^{I_\gamma} \not\equiv 1.
\]
\end{proposition}
\begin{proof} We fix $h_0\in \crc$ and a loop $\gamma\subset E_{h_0}$. The function $t$ from Lemma \ref{lem:tholo} together with $H$ provide the coordinates $(t,h)$ on the sufficiently small neighborhood of $\gamma$. We check by direct calculation that the function 
\[
 g(t) = \Big(e^{-I_\gamma}-1\Big)^{-1} \int_t^{t+T_\gamma} b(s)\cdot e^{-\int_t^s A}\,\dr s
\]
is periodic on $\gamma$ i.e.\ $g(t+T_\gamma) = g(t)$. This function is holomorphic in $t$ and meromorphic in $h$. The possible poles in $h$ are generated by the factor $\big(e^{-I_\gamma}-1\big)^{-1}$. The polynomial $p$ as well as the Hamiltonian time 1-form are holomorphic on the considered neighborhood of $\gamma$.

Thus the restriction of the integrated form $(g\, p\,\dr t)|_{E_h}$ is holomorphic on the open subset of the complex curve $E_h$. So, the integral $\Psi_\gamma$ is invariant under continuous deformations of $\gamma$ inside the fibre $E_h$. 

The resulting function $\Psi_\gamma (h)$ is meromorphic due to the analytical properties of $g$.\\
\end{proof}


\subsection{Local description of the generalized Abelian integrals}

\newcommand{\trg}{\mathbb{T}}
\newcommand{\trel}{W}

Below we introduce a relationship between the generalized Abelian integrals and certain representation of the fundamental group of the level curves of the polynomial $H$. This description explains the monodromy of the generalized Abelian integrals.

In the following Proposition-Definition we define certain $\Ad$-invariant function on the group of upper-triangular $3\times 3$ matrices.
\begin{proposition}
\label{pr:psfun}
Let $\trg$ be a group of upper-triangular matrices of the form
\begin{equation}
  \label{gform}
\trg \colon = \left\{ \trel=\begin{pmatrix}
    e^{-I/2}& a &b\\
    0& e^{I/2}& c \\
    0& 0& e^{-I/2} \\
    \end{pmatrix}, \quad I,a,b,c\in\bbC \right\}.
\end{equation}
The function $\psi(\trel)$, defined by the relation
\begin{equation}
  \label{psirdef}
  \big(e^{-I}-1\big)^{-1}\, \big[\trel- e^{I/2}\, \id \big]\,\big[\trel - e^{-I/2}\, \id\big] = \psi(\trel) \, \left(\begin{smallmatrix}0&0&1\\ 0&0&0\\ 0&0&0\\ \end{smallmatrix}\right),\qquad \trel\in\trg,
\end{equation}
is $\Ad$-invariant on the subset
\[
\mathbb{S}\colon = \{\trel\in\trg:\ (\det \trel)^2 \neq 1\}\subset \trg,
\]
on which $\psi(\cdot)$ is well defined.
\end{proposition}
\begin{remark}
We have the following explicit formula for the function $\psi$
\[
\psi\left(\begin{smallmatrix}e^{-I/2}&a &b\\ 0&e^{I/2}&c\\ 0&0& e^{-I/2}\\ \end{smallmatrix}\right) = e^{-I/2} b+ (e^{-I}-1)^{-1} a c.
\]
\end{remark}
\begin{remark}
The function $\psi(\cdot)$, defined by (\ref{psirdef}), is the obstruction to diagonalize the matrix $\trel$. In fact, provided $\psi(\trel)$ is well defined (i.e.\ $\trel\in \mathbb{S}$), the matrix $\trel$ is diagonalizable iff $\psi(\trel)=0$.
\end{remark}
\begin{proof}[Proof of Proposition \ref{pr:psfun}] Let $\widetilde{\trel}=\Ad_V \trel$, where $\trel, V$ (and also $\widetilde{\trel}$) are the matrices of the form (\ref{gform}). Since the adjoint transformation does not change eigenvalues, the diagonals in $\trel$ and $\widetilde{\trel}$ are the same.

The expression on the left hand side of the identity (\ref{psirdef}) is a quadratic polynomial $P_I(\cdot)$ evaluated on the matrix $\trel$. It satisfies $P_I(\widetilde{\trel}) = \Ad_V \big(P_I(\trel)\big)$. 

On the other hand, the polynomial $P_I(\trel)$ is proportional to the nilpotent matrix $\left(\begin{smallmatrix}0&0&1\\ 0&0&0\\ 0&0&0\\ \end{smallmatrix}\right)$ which commutes with any matrix from $\trg$. Thus
\[
P_I(\widetilde{\trel}) = \Ad_V \big(P_I(\trel)\big) = P_I(\trel).
\]
\end{proof}

\paragraph{The upper-triangular representation.}

Now we give a description of local properties of the generalized Abelian integral, under the general scheme presented in Proposition \ref{pr:ci}. We fix the point $h_0\in\crc$, corresponding complex level curve $E_0=H^{-1}(h_0)$ with chosen base point $*\in E_0$. The open neighborhood $U\subset \crc$ of $h_0$ is contractible.

We consider the following integrals along a loop $\gamma$ from the fundamental group $\pi_1(E_0,*)$:
\begin{equation}
\label{thtpdef}
  \begin{split}
    \tht_\gamma^+ &\colon =  e^{-I_\gamma/2}\,\int_0^{T_\gamma}\dr t\; p(t)\cdot e^{\int_0^t A},  \\
    \tht_\gamma^- &\colon =  e^{I_\gamma/2}\,\int_0^{T_\gamma}\dr t\; b(t)\cdot e^{-\int_0^t A}, \\
    \phi_\gamma   &\colon =  e^{-I_\gamma/2}\,\int_0^{T_\gamma}\dr t\int_0^t\dr s\; p(t)\cdot b(s)\cdot  e^{-\int_t^s A},
  \end{split}
\end{equation}
where the definition (\ref{igamdef}) of $I_\gamma$ and $T_\gamma$ is used.
\begin{remark}
It is worth to point out that the integrals $\tht_\gamma^\pm$ and $\phi_\gamma$ are well defined on the fundamental group $\pi_1(E_0,*)$ but, since they \emph{are not} $\Ad$-invariant, they do not give separately any contour integral. 
\end{remark}
\begin{theorem}
\label{th:rep}
The correspondence
\[
\pi_1(E_h,*)\ni \gamma \overset{\rho}{\longmapsto}
\begin{pmatrix}
  e^{-I_\gamma/2}& \tht_\gamma^- & \phi_\gamma\\
  0& e^{I_\gamma/2}& \tht^+_\gamma \\
  0& 0& e^{-I_\gamma/2} \\
\end{pmatrix}\in \trg
\]
gives rise to a \emph{group homomorphism}.
\begin{enumerate}
\item The subset 
  \begin{equation}
    \label{sdef}
    S\colon = \{\gamma \in \pi_1(E_0,*): \quad e^{I_\gamma} \not\equiv 1\}
  \end{equation}
is $\Ad$-invariant and the image of $S$ under the canonical projection onto $\pi_1(E_0)$ is monodromy invariant.
\item The contour integral $\Psi$, defined in (\ref{Psidef}), defines the $\Ad$-invariant function on $S$, which coincides with 
\[
\psi\circ\rho.
\]
\end{enumerate}
\end{theorem}
\begin{remark}
Note that the subset $S$ defined by (\ref{sdef}) is \emph{not} a group; the contractible loop does not satisfy the condition (\ref{sdef}).
\end{remark}
\begin{remark}
In the applications of generalized Abelian integrals to the analysis of the limit cycles of multi-dimensional dynamical systems (\ref{mdsys}) the condition $e^{I_\gamma} \not\equiv 1$ is guaranteed by the normal hyperbolicity -- see \cite{bozoell,bozo2}.
\end{remark}

\begin{proof}[Proof of Theorem \ref{th:rep}]
To prove that the map $\rho$ is a group homomorphism we have to verify the following integral identities
\begin{align}
  \tht^\pm_{\gamma\cdot\delta} &= e^{\pm I_{\gamma}/2}\,\tht^\pm_{\delta} + e^{\mp I_{\delta}/2}\,\tht^\pm_{\gamma}, \label{thtid}\\
  \phi_{\gamma\cdot\delta} &= e^{- I_{\gamma}/2}\,\phi_\delta + e^{- I_{\delta}/2}\,\phi_\gamma + \tht^-_\gamma\,\tht^+_\delta. \label{phid}
\end{align}
We calculate the integrals directly by the definition (\ref{thtpdef}), using the pa\-ra\-me\-tri\-zation of loops by the Hamiltonian time (recall that the period $\int_\gamma \dr t$ is denoted by $T_\gamma$):
\begin{align*}
  \tht^+_{\gamma\cdot\delta} &= e^{-(I_{\gamma}+I_{\delta})/2}\int_0^{T_{\gamma}+T_{\delta}}\dr t\; \Big(p(t)\cdot e^{\int_0^t A}\Big)=\\
&= \left[e^{-(I_{\gamma}+I_{\delta})/2}\int_0^{T_{\gamma}}\dr t + e^{-(I_{\gamma}+I_{\delta})/2} e^{I_{\gamma}}\int_0^{T_{\delta}}\dr t\right] \Big(p\, e^{\int A} \Big)\\
&= e^{I_{\gamma}/2}\,\tht^+_{\delta} + e^{- I_{\delta}/2}\,\tht^+_{\gamma}.
\end{align*}
This proves (\ref{thtid}) for the $+$ upper index. The calculations in the $-$ case are analogous.

To derive the composition rule (\ref{phid}) we calculate:
\begin{align*}
\phi_{\gamma\cdot\delta} &=  e^{-(I_\gamma+I_\delta)/2}\, \int_0^{T_\gamma+T_\delta}\dr t \int_0^{t}\dr s\; \left(p(t)\cdot b(s)\cdot e^{-\int_t^s A}\right)\\
 &= e^{-(I_\gamma+I_\delta)/2}\, \left[\int_0^{T_\gamma}\int_0^{t} + \int_{T_\gamma}^{T_\gamma + T_\delta}\int_{T_\gamma}^{t} + \int_{T_\gamma}^{T_\gamma + T_\delta}\int_0^{T_\gamma} \right] \Big(p\, b\, e^{-\int A} \Big)\\
&= e^{- I_{\delta}/2}\,\phi_\gamma + e^{- I_{\gamma}/2}\,\phi_\delta + e^{-(I_\gamma+I_\delta)/2}\,\int_{T_\gamma}^{T_\gamma + T_\delta}\dr t\int_0^{T_\gamma}\dr s \Big(p\, b\, e^{-\int A} \Big).
\end{align*}
Performing the coordinate change $t\mapsto t-T_\gamma$ and using the identity
\[
\int_{T_\gamma+t}^s A = \int_0^s A - I_\gamma - \int_0^t A
\]
we get
\begin{multline*}
e^{-(I_\gamma+I_\delta)/2}\,\int_{T_\gamma}^{T_\gamma + T_\delta}\dr t\int_0^{T_\gamma}\dr s \left(p(t)\cdot b(s)\cdot e^{-\int_t^s A}\right) = \\
= e^{(I_\gamma-I_\delta)/2}\,\int_0^{T_\delta}\dr t\int_0^{T_\gamma}\dr s \left(p(t)\,e^{\int_0^t A}\right)\cdot \left(b(s)\, e^{-\int_0^s A}\right) = \tht^-_\gamma\, \tht^+_\delta,
\end{multline*}
what finishes the proof of (\ref{phid}).
\medskip

To prove that the subset $S$ defined by (\ref{sdef}) is $\Ad$-invariant and monodromy invariant, we notice that the function $e^{I_\gamma}$ is the holomorphic contour integral whose domain of definition is whole $\pi_1(E)$. The condition $e^{I_\gamma} \not\equiv 1$ excludes some connected components of the covering $\pi_1(E)$ and so the remaining subset is monodromy invariant (see Remark \ref{rk:mon}). The subset $S$ is the pre-image under the canonical map (\ref{canpr}) of these components of $\pi_1(E)$. Thus $S$ is $\Ad$-invariant.

Now it remains to prove that the contour integral $\Psi$ defined by (\ref{Psidef}) coincides with $\psi\circ\rho$. We calculate
\begin{align*}
\Psi_\gamma =& \left(e^{-I_\gamma}-1\right)^{-1} \int_0^{T_\gamma}\dr t \int_t^{t+T_\gamma}\dr s\; \left(p(t)\cdot b(s)\cdot e^{-\int_t^s A}\right)\\
 =& \left(e^{-I_\gamma}-1\right)^{-1} \cdot\left[\int_0^{T_\gamma}\int_t^{T_\gamma} + e^{-I_\gamma}\int_0^{T_\gamma}\int_0^t \right] \Big(p\, b\, e^{-\int A} \Big) \\
=&\left[\int_0^{T_\gamma}\int_0^t + \left(e^{-I_\gamma}-1\right)^{-1} \int_0^{T_\gamma}\int_0^{T_\gamma}\right] \Big(p\, b\, e^{-\int A} \Big) =\\
\intertext{(by definition of the integrals $\tht^\pm_\gamma$ and $\phi_\gamma$)}
=& e^{I_\gamma/2}\, \phi_\gamma + \frac{\tht^+_\gamma\cdot\tht^-_\gamma}{e^{-I_\gamma} -1}.
\end{align*}
The latter expression coincides with $\psi\circ\rho$.

This finishes the proof of Theorem \ref{th:rep}.\\
\end{proof}

\subsection{The monodromy identities}
\label{sec:cor}

As a profit from the Theorem \ref{th:rep} one can easily get the monodromy action on the generalized Abelian integral. Quite complicated formulas (see \cite{bozoell,bozo2}) are consequences of some identity among certain functions on the group $\trg$ (see Lemma \ref{lem:mid} below). In addition to the generalized Abelian integral $\Psi$ we shall also need another contour integral
\[
\widetilde{\Psi}_\gamma = (e^{-I_\gamma} - 1)\cdot \Psi_\gamma
\]
and the corresponding function on the group $\trg$ defined by the relation (compare \ref{psirdef})
\[
 \big[\trel- e^{I/2}\, \id\big]\cdot \big[\trel - e^{-I/2}\, \id \big] = \widetilde{\psi}(\trel) \, 
\left(\begin{smallmatrix}0&0&1\\ 0&0&0\\ 0&0&0\\ \end{smallmatrix}\right).
\]
\begin{lemma}
\label{lem:mid}
Let $\trel_1,\trel_2\in\trg$ be such that $\psi (\trel_1)$, $\psi (\trel_2)$, $\psi (\trel_1\cdot \trel_2)$ are well defined. Then the following identity holds
\[
\psi(\trel_1\cdot \trel_2) = \psi(\trel_1)+ \psi(\trel_2) + \frac{e^{-(I_1+I_2)}}{(e^{-I_1} - 1)(e^{-I_2} - 1)(e^{-(I_1+I_2)} - 1)}\cdot \widetilde{\psi}([\trel_1,\trel_2]),
\]
where $[\trel_1,\trel_2] = \trel_1 \trel_2 \trel_1^{-1} \trel_2^{-1}$ is the commutant.

The analogous identities for the contour integrals $\Psi,\widetilde{\Psi}$ hold.
\end{lemma}
\begin{proof}
The above formula is a result of direct calculations.\\
\end{proof}

Now, knowing the monodromy of the loop, one can easily calculate the monodromy of generalized Abelian integrals. As an example we calculate the monodromy (and so singularity) of generalized Abelian integral in the elliptic case.

\paragraph{The elliptic case}

Let us consider the elliptic Hamiltonian $H=y^2-x^3+3x$. We assume that $A=\mathrm{const}$ and the polynomials $p,b$ are arbitrary.  The only critical values of $H$ are $\pm 2$, so $\crc=\bbC\setminus \{\pm 2\}$. The level curve $E_h$ for $h\in\crc$ is the elliptic curve (torus) with one point removed. As $h\in(-2,2)$, the curves $E_h$ contain the real oval $\gamma$ which we choose as a base point distinguishing a component of the covering $\pi_1(E)$ and so defining a multi-valued function $\Psi_\gamma$, meromorphic on the universal covering $\widetilde{\crc}$. The monodromy of the ``first branch'' of $\Psi_\gamma$ around $-2$ is trivial; the contour $\gamma$ squeezes to a critical point $(-1,0)$ as $h\to-2$. We determine the monodromy of $\Psi_\gamma$ around $h=+2$.

Let $\delta$ be a loop vanishing at the critical point $(1,0)$ as $h\to +2$. It is known (see \cite{zolbook,bozoell}) that the monodromies of contours are the following 
\begin{equation}
\label{monid}
\begin{split}
\mon_2 (\gamma) &= \gamma\cdot \delta, \\
\mon_2 (\delta) &= \delta, 
\end{split}
\end{equation}
what implies also
\begin{equation}
\label{monkid}
\mon_2 \big([\gamma,\delta]\big) = [\gamma,\delta].
\end{equation}

Using Theorem \ref{th:rep}, Lemma \ref{lem:mid} and the above monodromy relations we get
\begin{equation}
  \label{monell}
  \mon_2 \Psi_\gamma = \Psi_\gamma + \Psi_\delta - \frac{e^{-I_\delta}}{(e^{-I_\delta} - 1)^2}\left(\frac 1{e^{-(I_\gamma+I_\delta)} - 1} - \frac 1{e^{-I_\gamma} - 1}\right)\cdot  \widetilde{\Psi}_{[\gamma,\delta]}.
\end{equation}
We check, using the above formula, that the following auxiliary function
\[
\xi = \Psi_\gamma - \frac{\Psi_\delta}{2\pi i} \log(h-2) + \frac{e^{-I_\delta}\, \widetilde{\Psi}_{[\gamma,\delta]}}{(e^{-I_\delta} - 1)^2}\cdot \frac 1{e^{-I_\gamma} - 1}
\]
is monodromy invariant. Moreover, using the regular growth of functions involved in the definition of $\xi$, we deduce that $\xi$ is a germ of meromorphic function. Thus, the function $\Psi_\gamma$ has singularity of the following form
\[
\Psi_\gamma = \frac{\Psi_\delta}{2\pi i} \log(h-2) - \frac{e^{-I_\delta}\, \widetilde{\Psi}_{[\gamma,\delta]}}{(e^{-I_\delta} - 1)^2}\cdot \frac 1{e^{-I_\gamma} - 1} + \xi(h).
\]

The above results coincide with these from \cite{bozoell}, but now are achiev\-ed with much less calculations. 


Furthermore, one can prove the following
\begin{proposition}
\label{pr:ellneq}
The monodromy operator $\mon_2$, acting on the generalized Abelian integral $\Psi_\gamma$, generates infinite dimensional functional space. So, the function $\Psi_\gamma$ does not satisfy any linear differential equation of finite order.
\end{proposition}
\begin{remark}
The same should be true for the generalized Abelian integral of any Hamiltonian degree at least 3.
\end{remark}
\begin{remark}
The above proposition contrasts to either the Abelian integral or a generating function does (see Theorem \ref{th:de} and \cite{gi}).
\end{remark}

\noindent
\begin{proof}[Proof of Proposition \ref{pr:ellneq}]
Applying $(\mon_2 - \mathrm{Id})$ to the monodromy relation (\ref{monell}) we get
\[
\left(\mon_2 - \mathrm{Id}\right)^2 \Psi_\gamma =  -\frac{e^{-I_\delta}\, \widetilde{\Psi}_{[\gamma,\delta]}}{(e^{-I_\delta} - 1)^2} \left(\mon_2 - \mathrm{Id}\right)^2 \left(\frac 1{e^{-I_\gamma} - 1}\right),
\]
where the factor $ \frac{e^{-I_\delta}\, \widetilde{\Psi}_{[\gamma,\delta]}}{(e^{-I_\delta} - 1)^2}$ is a germ of meromorphic function, thanks to the monodromy invariance of $\delta$ and $[\gamma,\delta]$ -- see (\ref{monid}, \ref{monkid}). Since monodromy operator $\mon_2$ applied on $\left(\frac 1{e^{-I_\gamma} - 1}\right)$ generates an infinite dimensional vector space, the same is true for $\Psi_\gamma$.\\
\end{proof}


\begin{thebibliography}{Span}
\bibitem[AI]{aril}  V. I. Arnold and Yu. S. Il'yashenko, \emph{Ordinary differential equations, in: ``Ordinary Differential Equations and Smooth Dynamical Systems''}, Springer--Verlag, New York, (1997), pp. 1--148; (Russian: Fundamental Directions, v. \textbf{1}, {VINITI, Moscow}, 1985, pp. 1--146).
\bibitem[BZ1]{bozoell}
M. Bobie\'{n}ski and H. \.{Z}o\l \k{a}dek \emph{Limit cycles for multidimensional vector field. The elliptic case}, {J. Dynam. Control Systems}, {\bf 9}, {(2003)}, No 2, {265--310}.
\bibitem[BZ2]{bozo2}
M. Bobie\'{n}ski and H. \.{Z}o\l \k{a}dek \emph{Limit cycles of three dimensional polynomial vector fields}, Nonlinearity {\bf 18} 
(2005), No 1, 175--209.
\bibitem[H]{hall} M. Hall, \emph{The theory of groups}, New York, 1959.
\bibitem[FFG]{ffg}
{ A. T. Fomenko and D. B. Fuchs and V. L. Gutenmacher}, \emph{Homotopic topology}, {Akad\'emiai Kiad\'o (Publishing House of the Hungarian Academy of Sciences)}, {Budapest}, {1986}.
\bibitem[F1]{jpfr}
{J.-P. Fran{\c{c}}oise}, \emph{Successive derivatives of a first return map, application to the study of quadratic vector fields}, {Ergodic Theory Dynam. Systems}, {\bf 16}, {(1996)}, No 1,87--96.
\bibitem[F2]{jpfr2}
{J.-P. Fran{\c{c}}oise}, \emph{The successive derivatives of the period function of a plane vector field}, {J. Differential Equations}, {\bf 146}, {(1998)}, No 2, 320--335.
\bibitem[G]{newgav} L. Gavrilov, \emph{Higher order Poincar\'e-Pontryagin functions and iterated path integrals}, preprint, Laboratoire de Mathematiques Emile Picard (CNRS UMR 5580), 2004.
\bibitem[GI]{gi} L. Gavrilov and I. D. Iliev, \emph{The displacement map associated to polynomial unfoldings of planar Hamiltonian vector fields}, preprint, Laboratoire de Mathematiques Emile Picard (CNRS UMR 5580), 2003.
\bibitem[M]{mel} V. K. Melnikov, \emph{On the stability of a center for time-periodic perturbations}, Trans. Moscow Math. Soc. {\bf 12} (1963), 1--57.
\bibitem[LZ]{lezo}  P. Leszczy\'{n}ski, H. \.{Z}o\l \k{a}dek, \emph{Limit cycles appearing after perturbation of certain multi-dimensional vector fields}, J. Dynam. Diff. Equat. \textbf{13} (2001), No 4, 689--709.
\bibitem[P]{pont} L. S. Pontryagin, \emph{On dynamical systems close to Hamiltonian systems}, in: ``Selected Works'', v. 1, Gordon \& Breach, New York, 1986; [Russian: Zh. Ekper. Teoret. Fiziki 4 (1934), 234--238].
\bibitem[Z1]{zolq} H. \.Zo{\l}\k{a}dek, \emph{Quadratic systems with center and their perturbations}, {J. Differential Equations} \textbf{109} (1994), No 4, 223--273.
\bibitem[Z2]{zolbook} H. \.Zo{\l}\k{a}dek, \emph{The Monodromy Group}, Monografie Matematyczne, Birkh\"auser (submitted).
\end{thebibliography}
\end{document}